\documentclass[aps,apl,groupedaddress, showkeys,showpacs,twocolumn]{revtex4-1}
\usepackage[utf8]{inputenc}
\usepackage{amsmath,amssymb}
\usepackage{graphicx,float}
\makeatletter
\usepackage{dsfont}
\usepackage{multirow}
\title{GGT kinetiks}
\makeatother

\usepackage[bottom]{footmisc}

\newcommand{\iu}{\mathrm{i}\mkern1mu}

\usepackage{color}

\let\oldtop\top
\renewcommand{\top}{\oldtop\!}

\begin{document}

\title{Finite connected components in infinite directed and multiplex networks with arbitrary 
 degree distributions}
\author{Ivan Kryven} 
\email{i.kryven@uva.nl}
\affiliation{University of Amsterdam, PO box  94214, 1090 GE, Amsterdam, The Netherlands}

\begin{abstract}
This work presents exact expressions for size distributions of weak/multilayer connected components in two generalisations of the configuration model: networks with directed edges and multiplex networks with arbitrary number of layers. 
The expressions are computable in a polynomial time, and, under some restrictions, are tractable from the asymptotic theory point of view.
If first partial moments of the degree distribution are finite, the size distribution for two-layer connected components in multiplex networks exhibits exponent $-\frac{3}{2}$ in the critical regime, whereas  the size distribution of weakly connected components in directed networks exhibits two critical exponents, $-\frac{1}{2}$ and $-\frac{3}{2}$.

\end{abstract}
	

\maketitle

\newpage
\section{Introduction}

Many real world networks are well conceptualised when reduced to a graph, that is a set of nodes that are connected with edges or links. This representation helps to uncover often a non-trivial role of the topology in the functioning of complex networks \cite{kivela2014,nicosia2013,Kleineberg2016,cardillo2013}.
From a probabilistic perspective, many interesting network properties are well defined even when the total number of nodes approaches infinity.
For instance, the degree distribution is a univariate function of a discrete argument that denotes the probability for a randomly chosen node to have a specific number of adjacent edges \cite{newman2002}. The notion of degree distribution is easy to adapt to various generalisations of simple graphs. When different types of edges are present, or if edges are non-symmetrical (directed network), the degree distribution denotes the joint probability for a randomly sampled node to have specific numbers of edges of each type \cite{kivela2014}.

Just as a degree distribution is attributed to a single instance of a network, one may reverse this association, and talk about a class of networks that all match a given degree distribution. 
The class of such networks is known as the configuration model or generalised random graph \cite{bender1978,molloy1995,Newman2001,Newman2007,shi2009}.
In the configuration model, the connections between nodes are assigned at random with the only constraint that the degree distribution has to be preserved. This concept can naturally be extended to directed graphs, in which case the degree distribution is bivariate: counting incoming and outgoing edges \cite{Newman2001, kryven2016a}, or to multiplex networks where many types of edges exist, and thus, the degree distribution is multivariate \cite{bianconi2015,kivela2014,azimi2014,cellai2013,bianconi2013}.

A connected component is a set of nodes in which each node is connected to all other nodes with a path of finite or infinite length. Different notions of a path give rise to distinct definitions of connected components. Namely, if directed edges are present: in-, out-, weak, and strong components are distinguished \cite{Newman2001}. As in multiplex networks, one may speak of a connected component that is solely contained within a single layer or a two-layer component having edges in both layers \cite{bianconi2015,baxter2014,Domenico2013}.
Even under the assumption of the thermodynamic limit, when the total number of nodes approaches infinity, the infinite network may contain connected components of finite size $n>1$.
And so there are two key features that characterise sizes of connected components in configuration models: the size distribution of finite components and the size of the giant component. The size distribution is usually defined as the probability that a randomly sampled node belongs to a component of a specific size, while the size of the giant component is the probability that a randomly sampled node belongs to a component of size that scales linearly with the size of the whole system \cite{Newman2001}. 

A considerable progress has been made in recovering both the size distribution and the size of the giant component that are associated with an arbitrary degree distribution in undirected, single-layer configuration networks.
Molloy and Reed \cite{molloy1995} proposed a simple criterion to test the existence of the giant component.
In Ref.~\cite{Newman2001}, Newman \emph{et al.} narrowed the problem of finding the size distribution down to a numerical solution of an implicit functional equation, that is followed by the generating function inversion. Somewhat later, a few cases have been resolved analytically \cite{Newman2007}, and recently, the formal solution for size distribution of connected components in undirected networks has been found by means of the Lagrange inversion \cite{kryven2017a,kryven2017b}. Such a solution permits fast computation of exact numerical values and allows simple asymptotic analysis. 

A smaller amount of results, however, is available for \emph{directed} and \emph{multiplex} configuration models. In these cases the aforementioned functional equation remains the main bottleneck and is typically addressed numerically with the only exception of percolation studies. Some percolation criteria were obtained analytically both in directed networks: in-/out-percolation \cite{Newman2001}, weak percolation \cite{kryven2016a}, and in multiplex networks: $k$-core percolation \cite{azimi2014}, weak percolation \cite{baxter2014}, strong mutually connected component \cite{bianconi2015} and giant connected component \cite{hackett2016}. Up to date, little results are available on the \emph{size distribution} of finite connected components in these configuration models.

The present paper applies the Lagrange inversion principle to find exact expressions for size distributions of connected components in two generalisations of the configuration model: directed configuration networks and multiplex configuration networks. Firstly, a brief review of the Good's multivariate generalisation of the Lagrange inversion formula is given. Then, the size distributions for in-, out-, and weak components in directed configuration networks are formulated in terms of convolution powers of the degree distribution. These results are complemented by a detailed asymptotic analysis that reveals existence of two distinct critical exponents. In the next section, the general case of weak multi-layer connected components (i.e. components that include edges from an arbitrary layer) is considered. 
A formal expression for the size distribution is constructed and the asymptotic analysis is provided for two-layer multiplex networks. Furthermore, the relation between these results and the existence of a two-layer giant component is studied by means of perturbation analysis within the critical window. Finally, the results for directed and multiplex networks are illustrated with a few examples in the last section.

\section{Lagrange series inversion}
Suppose $R(x),A(x),F(x)$ are formal power series in $x.$  Then, according to the Lagrange inversion formula \cite{bergeron1998}, implicit functional equation 
\begin{equation}
\label{eq:SimpleImplicit}
A(x)=x R[A(x)]
\end{equation}
 has a unique solution $A(x)$. Instead of an expression for $A(x),$ the Lagrange inversion formula recovers a discrete function that is generated by $A(x).$ In fact, the equation yields a slightly more general result: for an arbitrary formal power series $F(x),$ the coefficients of power series $F[A(x)]$ at $x^n$ read as,
\begin{equation}
\label{eq:LagrangeInversion_1d}
[x^n]F[A(x)] = \frac{1}{n}[t^{n-1}]F'(t)R^n(t),\; n>0.
\end{equation}
Here $[t^{n-1}]$ refers to the coefficient at $t^{n-1}$ of the corresponding power series. 
In the context of configuration models, Eq.~\eqref{eq:LagrangeInversion_1d} proved to be useful when deriving a formal expression for the size distribution of connected components in undirected networks \cite{kryven2017a}.

The Lagrange inversion, was generalised to the case of multivariate series by Good \cite{good1960}. Following the original notation  from \cite{bergeron1998}, the Lagrange-Good theorem in $d$ dimensions reads: let $\mathbf{R}(\mathbf{x}) = [R_1(\mathbf{x}),R_2(\mathbf{x}),\dots,R_d(\mathbf{x})]$ be a vector of formal power series in variables $\mathbf{x}=(x_1,x_2,\dots,x_d),$ and let $\mathbf{A}(\mathbf{x})$ be a vector of formal power series satisfying 
\begin{equation}
\label{eq:Implicit_ND}
\mathbf A_i(x_1,\dots,x_d)=x_i R_i(\mathbf A_1,\dots,\mathbf A_d), \; i=1,\dots,d,
\end{equation}
 then for any formal power series $ F(\mathbf x)$, 
\begin{equation}
 \label{eq:LagrangeInversion_nd}
 [\mathbf x^\mathbf n] F[\mathbf A(\mathbf x)]=[\mathbf t^\mathbf n]  F(\mathbf t)\text{det}[ K (\mathbf t)]\mathbf R^\mathbf n(\mathbf t),\;\mathbf n\in \mathbb{N}^d,
\end{equation}
where $K(\mathbf t)$ is a matrix from $\mathbb{R}^{d\times d},$
\begin{equation}
\label{eq:NDdet}
\tag{\ref{eq:LagrangeInversion_nd}$a$}
K(\mathbf t)_{i,j}=\delta_{i,j} - \frac{t_i}{\mathbf R_i(\mathbf t)} \frac{\partial \mathbf R_i}{\partial t_j}(\mathbf t),\; i,j=1,\dots,d,
\end{equation}
and
$\mathbf t= (t_1,\dots,t_d),$ $\mathbf n= (n_1,\dots,n_d),$  $\mathbf x^\mathbf n = [x_1^{n_1},\dots,x_d^{n_d}],$ $\mathbf x(\mathbf y) = [x_1(\mathbf y),\dots x_d(\mathbf y)].$ Analogously to the one-dimensional case \eqref{eq:LagrangeInversion_1d}, the operator $[\mathbf x^\mathbf n ]$ refers to the coefficient at $x_1^{n_1},\dots,x_d^{n_d}.$ In the case when $d=1$, Eq.~\eqref{eq:LagrangeInversion_nd} simplifies to the Lagrange equation \eqref{eq:LagrangeInversion_1d}. Although the original formulation of the Langrange-Good equation \eqref{eq:LagrangeInversion_nd} does involve an inversion of a generating function (GF), the only reason the inversion is used is to perform the convolution. Where convenient, we will  exploit this fact and write \eqref{eq:LagrangeInversion_1d} without any reference to GFs at all by utilising the convolution power notation: $f(\mathbf k)^{*n}=f(\mathbf k)^{*n-1} * f(\mathbf k),\; f(\mathbf k)^{*0}:=\delta(\mathbf k)$, where
the multidimensional convolution is defined as $d(\mathbf n)=f(\mathbf k)*g(\mathbf k ),$
\begin{equation}\label{eq:convnd}d(\mathbf n) = \sum\limits_{\mathbf j + \mathbf k=\mathbf n}f(\mathbf j)g(\mathbf k)=[t^{\boldsymbol k}] F(x)G(x).\end{equation}
Here, $\mathbf i,\mathbf j,\mathbf k,\mathbf n$ are $d-$dimensional vectors. The sum in Eq.~\eqref{eq:convnd} runs over all partitions of vector $\mathbf n$ into two summands $\mathbf j,\,\mathbf k,$ such that $$j_i+k_i=n_i,\;0\leq j_i,k_i \leq n_i, \; i=1,\dots,d.$$

In practice, numerical values of the convolution can be conveniently obtained with Fast Fourier Transform (FFT).
We will see now how the inversion equations \eqref{eq:LagrangeInversion_1d} and \eqref{eq:LagrangeInversion_nd} can be applied to find the size distributions for connected components in directed and multiplex networks that are defined by their degree distributions. 

\section{Directed networks}\label{eq:SectionDirrected}
In a directed network, bivariate degree distribution $0 \leq u(k,l)\leq1$ denotes probability of choosing a node with $k\geq0$ incoming edges and $l\geq0$ outgoing edges uniformly at random. Partial moments of this distribution are given by 
\begin{equation}\label{eq:mom}
\mu_{ij}=\sum\limits_{k,l=0}^\infty k^il^j u(k,l).
\end{equation}
Since $u(k,l)$ is normalised, $\mu_{00}=1,$ and since the expected numbers for incoming and outgoing edges must coincide, $\mu_{10}=\mu_{01} =\mu.$
Directed degree distribution $u(k,l)$ has two corresponding excess distributions: $u_{\text{in}}(k,l) =\frac{k+1}{\mu} u(k+1,l)$ and $u_{\text{out}}(k,l) =\frac{l+1}{\mu} u(k,l+1)$. 
Throughout this section, the capital letters are used to denote the corresponding bivariate GFs: $U(x,y),U_{\text{in}}(x,y),U_{\text{out}}(x,y).$ 
Four types of connected components are distinguished in directed configuration models: in-components, out-components, weak component, and strong component (the latter always has an infinite size in the thermodynamic limit \cite{Newman2001}).

\subsection{Sizes of in- and out-components}
 The size distributions for both, in-components $h_{\text{in}}(n),$ as generated by $H_{\text{in}}(x)$, and out-components $h_{\text{out}}(n),$ as generated by $H_{\text{out}}(x),$ can be found by solving the following systems of functional equations \cite{Newman2001}:
\begin{equation}
	\begin{aligned}
		H_{\text{out}}(x) =& x U\big[ \tilde H_{\text{out}}(x) ,1 \big],\\
		\tilde H_{\text{out}}(x) =&  xU_{\text{out}}[  \tilde H_{\text{out}}(x),1 ]
	\end{aligned}
\end{equation}
and
\begin{equation}
	\begin{aligned}
		H_{\text{in}}(x) =& x U\big[1,\tilde H_{\text{in}}(x)\big],\\
		\tilde H_{\text{in}}(x) =&  xU_{\text{in}}[1, \tilde H_{\text{in}}(x) ].
	\end{aligned}
\end{equation}
These equations are similar to those describing connected components in the undirected configuration network, and following a similar derivation to the one from Ref.~\cite{kryven2017a}, one immediately obtains formal solutions in terms of the convolution power of the degree distribution,
\begin{equation}
\label{eq:Lagrange1directed}
\begin{aligned}
&h_{\text{in}}(n)=\frac{\mu}{n-1} \tilde u_{\text{in}}^{*n}(n-2),  \;n>1;\\
&h_{\text{out}}(n)=\frac{\mu}{n-1} \tilde u_{\text{out}}^{*n}(n-2),  \;n>1;\\
&h_{\text{in}}(1)=h_{\text{out}}(1)=u(0,0).
\end{aligned}
\end{equation}
Here $\tilde u_{\text{in}}(k)=\sum\limits_{l=0}^{\infty} u_{\text{in}}(k,l)$ and $\tilde u_{\text{out}}(l)=\sum\limits_{k=0}^{\infty} u_{\text{out}}(k,l).$
\subsection{Weakly connected components}
The generating function for the size distribution of weak components $W(x),$ satisfies the following system of functional equations \cite{kryven2016a},
\begin{equation}\label{eq:W1}
W(x) = x U\Big[ W_{\text{out}}(x), W_{\text{in}}(x) \Big],
\end{equation}
\begin{equation*}
\tag{\ref{eq:W1}$a$}
\begin{aligned}
&W_{\text{out}}(x) = x U_{\text{out}}\Big[ W_{\text{out}}(x), W_{\text{in}}(x) \Big],\\
&W_{\text{in}}(x)  = x U_{\text{in}} \Big[ W_{\text{out}}(x), W_{\text{in}}(x) \Big].
\end{aligned}
\end{equation*}
To solve this system we apply the Lagrange-Good formalism \eqref{eq:Implicit_ND}. First, one should transform \eqref{eq:W1} to match the bi-variate version ($d=2$) of Eq.~\eqref{eq:Implicit_ND}. Consider three bi-variate formal power series, $A(x,y),A_1(x,y),A_2(x,y)$ that take their diagonals from correspondingly $\frac{1}{x}W(x),\;W_{\text{out}}(x),\;W_{\text{in}}(x),$ that is 
\begin{equation}
\label{eq:diag_def}
\begin{aligned}
&A(x,x)=\frac{1}{x}W(x),\\
&A_1(x,x)=W_{\text{out}}(x),\\
&A_2(x,x)=W_{\text{in}}(x), 
\end{aligned}
\end{equation}
for $|x|<1,\; x\in \mathbb{C}$.
Additionally, let $R_1(x,y):=U_{\text{out}}(x,y),\; R_2(x,y):=U_{\text{in}}(x,y).$ If  couple $A_1(x,y),A_2(x,y)$ satisfies condition \eqref{eq:Implicit_ND} for all values of $(x,y),$ then as being a partial case ($x=y$), the weaker condition  (\ref{eq:W1}a) is also satisfied. Furthermore, by assigning $F(x,y):= U(x,y)$ one obtains the expression for the coefficients of generating function $A(x,y)$: for $i,j\geq0$,
\begin{multline}
\label{eq:ajk0}
a(i,j)=[x^iy^j]A(x,y)=[x^iy^j]  U[ A_1(x,y), A_2(x,y) ] =\\ [t_1^i t_2^j]  U(t_1,t_2)\text{det}[ K (t_1,t_2)] U_{\text{out}}(t_1,t_2)^i U_{\text{in}}(t_1,t_2)^j,
\end{multline}
which when rewritten with the convolution power notation \eqref{eq:convnd}, become
\begin{equation}
\label{eq:ajk}
a(i,j)= u(k,l)* u_{\text{out}}(k,l)^{*i-1}* u_{\text{in}}(k,l)^{* j-1 } * d(k,l)\Big|_{\substack{k=i\\l=j}},
\end{equation}
where
\begin{multline}
\label{eq:Djk}
d(k,l)= [ u_{\text{out}}(k,l)- k u_{\text{out}}(k,l)] *  [ u_{\text{in}}(k,l) - l u_{\text{in}}(k,l) ] \\
- l u_{\text{out}}(k,l) * k u_{\text{in}}(k,l).  
\end{multline}
Here, $d(k,l)$ is chosen in such a way that it is generated by  $U_{\text{out}}(t_1,t_2)U_{\text{in}}(t_1,t_2)\det[K(t_1,t_2)],$ the product that appears in Eq.~\eqref{eq:ajk0}. For this reason the convolution powers in Eq.~\eqref{eq:ajk} are diminished by one: $i-1, \; j-1$.
Now, on one hand $w(n+1)$ is generated by $\frac{1}{x}W(x)=A(x,x),$  on the other $x^i y^j|_{y=x} = x^{i+j}$ and thus the sum of all $a(i,j)=[x^j y^j]A(x,y)$ such that $i+j=n+1$ yields the values of $w(n+1)$. Therefore, the final expression for the size distribution of  weak components is written out as a diagonal sum,
\begin{equation}
\label{eq:w2d}
w(n) =
\begin{cases}
\sum\limits_{i=0}^{n-1} a(i,n-i-1), & n>1;\\
 u(0,0), & n=1.
 \end{cases}
\end{equation}
From the computational perspective, the most efficient way to evaluate Eq.~\eqref{eq:ajk} numerically is to apply FFT algorithm to find the convolution powers. In this case, the computation of $w(n)$ requires $O(n^2 \log n)$ multiplicative operations. 

Besides being suitable for numerical computations, expressions \eqref{eq:Lagrange1directed} and \eqref{eq:w2d} can be further treated analytically to obtain the asymptotic behaviour of size distributions $w(n), h_{\text{in}}(n), h_{\text{out}}(n)$ in the large $n$ limit. 
That is we will search for such $w_\infty(n)$ (or correspondingly $h_{\text{in},\infty}(n)$ and $h_{\text{out},\infty}(n)$) that 
\begin{equation}\label{eq:def:winf}
\frac{w(n)}{w_\infty(n)}\to 1,\; n\to \infty. 
\end{equation}
In the context of asymptotic theory, we limit ourself to the case of finite first moments, $\mu_{ij}<\infty,\;i+j \leq 3.$ As will be shown further on, this assumption will allow us to utilise the standard central limit theorem and formulate the analytical expressions for the asymptotes as a function of solely the first partial moments of the degree distribution, $\mu_{ij},\; i+j\leq 3.$ 
To keep the derivation concise, we define shorthands for the vectors of expected values and covariance matrices of $u(k,l),$ $\frac{k}{\mu_{10}}u(k,l),$ and $\frac{l}{\mu_{01}}u(k,l):$
\resizebox{1\columnwidth}{!}{
  \begin{minipage}{1.28\columnwidth}
  \begin{equation}\label{eq:mu1mu2}
\begin{aligned}
\mu_0 =\begin{bmatrix} \mu_{10} \\ \mu_{01}
\end{bmatrix},
\;
&
\Sigma_0=\begin{bmatrix}
\mu_{20}- \mu_{10}^2              & \mu_{11} - \mu_{10} \mu_{01}  \\
\mu_{11} - \mu_{10} \mu_{01}   & \mu_{02}- \mu_{01}^2\\
\end{bmatrix};
\\
\mu_1 =\frac{1}{\mu_{10}}\begin{bmatrix} \mu_{20} \\ \mu_{11}
\end{bmatrix},
\;
&
\Sigma_1=\frac{1}{\mu^2_{10}}
\begin{bmatrix}
\mu_{30} \mu_{10}- \mu_{20}^2              & \mu_{21}\mu_{10} -\mu_{11} \mu_{20}  \\
\mu_{21} \mu_{10} -\mu_{11} \mu_{20}   & \mu_{12}\mu_{10}-\mu_{11}^2\\
\end{bmatrix};
\\
\mu_2 =\frac{1}{\mu_{01}}\begin{bmatrix} \mu_{11} \\ \mu_{02}
\end{bmatrix},\;
&
\Sigma_2=
\frac{1}{\mu^2_{01}}
\begin{bmatrix}
 \mu_{21} \mu_{01}-\mu_{11}^2  &  \mu_{12}\mu_{01} - \mu_{02} \mu_{11} \\
 \mu_{12}\mu_{01}-\mu_{02} \mu_{11} &  \mu_{03}\mu_{01}-\mu_{02}^2 \\
\end{bmatrix}.\\
\end{aligned}
\end{equation}
\end{minipage}
}\\

\noindent Note, that in directed networks $\mu_{10}=\mu_{01}=\mu.$

\subsection{Asymptotes for in- and out-components}
In the case of in- and out-components the asymptotic analysis coincides with the one performed in the case of undirected network and has been covered elsewhere, for instance, compare Eq.~\eqref{eq:Lagrange1directed} to Eq.~(8) in  Ref.~\cite{kryven2017a}. Taking this into the account, we can immediately proceed with expressions for the asymptotes: 
\begin{equation}\label{eq:h_in_inf}
\begin{aligned}
&h_{\text{in},\infty}(n) = C_{1,1}  e^{- C_{1,2} n}n^{-\frac{3}{2}},\\
& C_{1,1} = \frac{\mu^2}{ \sqrt{2 \pi (\mu \mu_{30}-\mu_{20}^2) }}, \;C_{1,2} = \frac{(\mu_{20}-2 \mu)^2}{2( \mu \mu_{30}- \mu_{20}^2)};\\
\end{aligned}
\end{equation}
\begin{equation}\label{eq:h_out_inf}
\begin{aligned}
&h_{\text{out},\infty}(n)= C_{2,1}  e^{-  C_{2,2}n}n^{-\frac{3}{2}},\\
& C_{2,1} =\frac{\mu^2}{ \sqrt{2 \pi (\mu \mu_{03}-\mu_{02}^2) }}, \;C_{2,2} = \frac{(\mu_{02}-2 \mu)^2}{2( \mu \mu_{03}- \mu_{02}^2)},
\end{aligned}
\end{equation}
and refer the reader to Ref.~\cite{kryven2017a} for the derivation. One can see that depending on the values of the moments, the asymptotes \eqref{eq:h_in_inf},\eqref{eq:h_out_inf} switch between exponential and algebraic decays. The algebraic asymptote exhibits slope $-\frac{3}{2}$, which implies that in this case the size distributions feature infinite expected values. 
According to Eqs. \eqref{eq:h_in_inf},\eqref{eq:h_out_inf}, the algebraic asymptote emerges when $\mu_{20}-2 \mu=0,$ for in-components, and $\mu_{02}-2 \mu=0,$ for out-components, both of which coincide with the critical point for the existence of the corresponding giant components \cite{kryven2016a}.

\subsection{Asymptote for weakly connected components}\label{sec:asdirect}
The asymptotic analysis for the size distribution of weak components is conceptually different from the previous case: unlike in Eq.~\eqref{eq:Lagrange1directed}, the expression for size distribution \eqref{eq:ajk}-\eqref{eq:w2d} contains the complete bivariate degree distribution and therefore cannot be treated analogously to the case of undirected networks.

We start by replacing the generating function appearing in the right hand side (RHS) of Eq.~\eqref{eq:ajk0} with 
a characteristic function by introducing a change of variables $t_1=e^{\iu \omega_1},\; t_2=e^{\iu\omega_2}$:
   \begin{multline}\label{eq:phia}
   \phi_a(\omega_1,\omega_2) =   U(e^{\iu \omega_1},e^{\iu \omega_2})\text{det}[ K (e^{\iu \omega_1},e^{\iu \omega_2})]\\
   \times U_{\text{out}}(e^{\iu \omega_1},e^{\iu \omega_2})^i U_{\text{in}}(e^{\iu \omega_1},e^{\iu \omega_2})^j.
   \end{multline}
Here, the complex unity is denoted with $\iu^2=-1$; it should not be confused with the parameters $i,j$.
By setting $\phi( \omega_1, \omega_2 ):= U( e^{\iu \omega_1}, e^{\iu \omega_2})$ and expanding $ U_{\text{in}},\, U_{\text{out}},$ and $K$ according to their definitions  one obtains 
\resizebox{1\columnwidth}{!}{
  \begin{minipage}{1.28\columnwidth}
\begin{equation}
\begin{aligned}
\label{eq:long_ch}
\phi_a(&\omega_1,\omega_2)= e^{-\iu (i \omega_2+ j\omega_1)} \phi(\omega_1, \omega_2)\\
\times  &\left( \left[-\frac{\iu}{\mu} \frac{\partial}{\partial \omega_1}\phi(\omega_1, \omega_2)\right]^j
\left[-\frac{\iu}{\mu} \frac{\partial}{\partial \omega_2}\phi(\omega_1, \omega_2)\right]^i\right.
 \\ 
&+\frac{1 }{j}  \iu
\frac{\partial }{\partial \omega_2}\left[- \frac{\iu}{\mu} \frac{\partial}{\partial \omega_1}\phi(\omega_1, \omega_2)\right]^j
\left[-\frac{\iu}{\mu} \frac{\partial}{\partial \omega_2}\phi(\omega_1, \omega_2)\right]^i\\
&+
\frac{1}{i}  
\left[-\frac{\iu}{\mu} \frac{\partial}{\partial \omega_1}\phi(\omega_1, \omega_2)\right]^j
\iu\frac{\partial }{\partial \omega_1}\left[-\frac{\iu}{\mu} \frac{\partial}{\partial \omega_2}\phi(\omega_1, \omega_2)\right]^i\\
&-
\frac{1}{i j} 
\frac{\partial }{\partial \omega_2}\left[- \frac{\iu}{\mu}\frac{\partial}{\partial \omega_1}\phi(\omega_1, \omega_2)\right]^j
\frac{\partial }{\partial \omega_1}\left[- \frac{\iu}{\mu}\frac{\partial}{\partial \omega_2}\phi(\omega_1, \omega_2)\right]^i\\
&+
\left.
\frac{1}{i j} 
\frac{\partial }{\partial \omega_1}\left[-  \frac{\iu}{\mu}\frac{\partial}{\partial \omega_1}\phi(\omega_1, \omega_2)\right]^j
\frac{\partial }{\partial \omega_2}\left[-\frac{\iu}{\mu} \frac{\partial}{\partial \omega_2}\phi(\omega_1, \omega_2)\right]^i \right).
\end{aligned}
\end{equation}
\end{minipage}
}\\

Having $\phi_a(\omega_1,\omega_2)$ in this format, allows us to apply the central limit theorem, that guarantees the pointwise convergence of the following limits:
\resizebox{1\columnwidth}{!}{
\begin{minipage}{1.28\columnwidth}
\begin{equation}\label{eq:limits}
\begin{aligned}
\lim_{j\to\infty}|\left(-\frac{\iu}{\mu} \frac{\partial}{\partial \omega_1}\phi(\omega_1, \omega_2)\right)^j-\phi_g( \omega, j\mu_1, j\Sigma_1)|=0,\\
\lim_{i\to\infty}|\left(- \frac{\iu}{\mu}\frac{\partial}{\partial \omega_2}\phi(\omega_1, \omega_2)\right)^i-\phi_g(\omega, i\mu_2, i\Sigma_2)|=0.\\
\end{aligned}
\end{equation}
\end{minipage}
}\\
Here $\phi_g(\omega,\mu,\Sigma)=e^{\iu\, \mu^\top\omega-\frac{1}{2} \omega^\top\Sigma\,\omega}, $ $\omega=(\omega_1,\omega_2)^\top$ denotes the characteristic function for the bivariate Gaussian-distributed random variable and $\mu_1,\mu_2,\Sigma_1,\Sigma_2$ are as defined in Eq.~\eqref{eq:mu1mu2}.
Now, after substituting the limiting functions from \eqref{eq:limits} into \eqref{eq:long_ch}, evaluating the partial derivatives, and using the symmetricity of matrices $\Sigma_1,\Sigma_2$, one obtains:
\resizebox{1\columnwidth}{!}{
\begin{minipage}{1.28\columnwidth}
\begin{equation}
\begin{aligned}\label{eq:ch1}
\phi_{a,\infty}(\omega_1,\omega_2)= & 
e^{\iu\, [ j\mu_1 +i\mu_2 - (j,i)+\mu_0]\omega-\frac{1}{2} \omega^\top(j\Sigma_1+i\Sigma_2)\,\omega}  \\ \times
 & \left( I( \mu_1-\mu_2) +\mu_1^\top D \mu_{2} \right.\\
 &+[I  (\Sigma_1-\Sigma_2)  + \mu_1^\top D  \Sigma_2 -\mu_2^\top D  \Sigma_1 ]\omega \\
 &-\left. \omega^\top \Sigma_2 D \Sigma_1\omega \right),
 \end{aligned}
 \end{equation}
 \end{minipage}
}\\
 where
 $$
 D=
\begin{bmatrix}
0 & \text{-}1\\
1 & \;0\\
\end{bmatrix},\; I=\begin{bmatrix}\;1\\\text{-1}\end{bmatrix}.$$
Note that Eq.~\eqref{eq:ch1} does not contain $\Sigma_0,$ which becomes negligible in the limit of large $i,j.$
After applying the inverse Fourier transform, \eqref{eq:ch1} becomes
\begin{equation}\label{eq:a_inf}
a_{\infty}(i,j) =C(x) \frac{ e^{ -\frac{1}{2} x^\top\Sigma^{-1} x }}{2 \pi  \sqrt{\det(\Sigma)}},\\
\end{equation}
where $x=(i,j)+(j,i) - ( j\mu_1^\top +i\mu_2^\top)+\mu_0$ and
$$
\begin{aligned}
C(x) = &C_0+C_1  \frac{x}{n} + x^\top C_2\frac{x}{n^2},\\
C_0    = &I( \mu_1-\mu_2) +\mu_1^\top D \mu_{2},\\
C_1    = &[I  (\Sigma_1-\Sigma_2)  + \mu_1^\top D  \Sigma_2 -\mu_2^\top D  \Sigma_1 ]\left(\frac{1}{n}\Sigma\right)^{-1}\!\!\!\!\!\!, \\
C_2    = &-\left(\frac{1}{n}\Sigma\right)^{-1} \Sigma_2 D \Sigma_1 \left(\frac{1}{n}\Sigma\right)^{-1}\!\!\!\!\!\!, \\
\Sigma=&j\Sigma_1+i\Sigma_2.
\end{aligned}
$$
By introducing new variables 
\begin{equation}\label{eq:new_variables}
n= i+j,\;z=\frac{i-j}{i+j},
\end{equation}
  one obtains:
$x=n (a  z+ b  + \mu_0/n )$ and
\begin{equation}\label{eq:Sigma}
\Sigma = n ( A z + B ),
\end{equation}
where
\begin{equation}
\label{eq:AB}
\begin{aligned}
&a = \frac{\mu_1-\mu_2}{2},\; &
b = 1- \frac{\mu_1 + \mu_2}{2},\\
&A=\frac{\Sigma_2- \Sigma_1}{2},\;&
B=\frac{\Sigma_1 + \Sigma_2}{2}.
\end{aligned}
\end{equation}
Under this change of variables, $\frac{1}{n}\Sigma=(Az+B)$  is independent of $n$ and, consequently, so are $C_0,C_1,C_2.$
Furthermore, the exponential function from \eqref{eq:a_inf} can be now rewritten as a univariate Gaussian function in $z,$
\begin{multline}
\frac{ e^{ -\frac{1}{2} x^\top(\Sigma)^{-1} x }}{2 \pi  \sqrt{\det(\Sigma)}}=\\
\frac{ e^{ -\frac{1}{2} n^2 (a  z+ b  + \mu_0/n )^\top[n (Az+B)]^{-1} (a  z+ b  + \frac{\mu_0}{n} ) }}{2 \pi  \sqrt{\det[n (Az+B)]}}=\\
\frac{ e^{ -\frac{1}{2}  (a  z+ b  + \mu_0/n )^\top(\frac{Az+B}{n})^{-1} (a  z+ b  + \frac{\mu_0}{n} ) }}{2 \pi \sqrt{\det[n (Az+B)]}}=\\
\frac{ e^{ -\frac{\left(z+ \frac{  \mu_0^\top S_n^{-1}a/n+ a^\top S_n^{-1}b }{a^\top S_n^{-1}a  } \right)^2  } {2(a^\top S_n^{-1}a)^{-1}  }}} {\sqrt{2 \pi 
(a^\top S_n^{-1}a)^{-1}
}}C_b(n,z)  =\\ C_b(n)\mathcal{N}(z,-\frac{  \mu_0^\top S_n^{-1}a/n+ a^\top S_n^{-1}b }{a^\top S_n^{-1}a  },a^\top S_n^{-1}a)
\end{multline}
where $S_n= \frac{A z + B }{ n },$
and 
$$C_b(n)= \frac{e^{ -\frac{1}{2}\left[   (b+\mu_0/n)^\top S_n^{-1} (b+\mu_0/n)  -\frac{ ( \mu_0^\top S_n^{-1}a/n+ a^\top S_n^{-1}b )^2}{a^\top S_n^{-1}a  } \right] }}{n^2 \sqrt{2 \pi   \det(S_n)\, a^\top S_n^{-1} a} },
$$
and $ |z|\leq 1.$
At the limit $n\to \infty,$ the variance of this Gaussian function vanishes as $O(n^{-1})$, and the expected value remains bounded. Indeed, for a fixed $z,$ such that $S_n^{-1}$ exists:
$$
\begin{aligned}
&a^\top S_n^{-1} a  =  O(n^{-1}),\\
&\frac{  \mu_0^\top S_n^{-1}a/n+ a^\top S_n^{-1}b }{a^\top S_n^{-1}a } =\frac{a^\top(Az+B)^{-1}b }{a^\top(Az+B)^{-1}a }+ O(n^{-1}),
\end{aligned}
$$
so that the Gaussian function itself tends to the Dirac delta function $\delta[z+\frac{a^\top(Az+B)^{-1}b }{a^\top(Az+B)^{-1}a } ].$
Recall, that according to \eqref{eq:w2d} the size distribution is defined as a sum of the diagonal elements  $$w_{\infty}(n+1)=\sum\limits_{i+j=n}a_{\infty}(i,j)=\sum\limits_{k=1}^n a_{\infty}(i,j)\Big|_{\substack{i+j=n,\;\;\;\;\;\;\\j=(i-j)k/n.}}$$  
This sum can be viewed as an estimator for an integral,
\resizebox{1\columnwidth}{!}{
\begin{minipage}{1.28\columnwidth}
 \begin{equation}
 \label{eq:integral}
 \begin{aligned}
 w_\infty(n+1) = \frac{1}{n}\int\limits_{-1}^1& \delta[z+\frac{a^\top(Az+B)^{-1}b }{a^\top(Az+B)^{-1}a } ]\\
& \times C_a[n (a  z+ b  + \frac{\mu_0}{n} )]C_b(n,z)\, \text{d}z,
 \end{aligned} 
 \end{equation}
 \end{minipage}
}\\
such  that $\lim\limits_{n\to\infty}|w(n) - w_\infty(n)|=0.$ The delta function under the integral is non-zero only at $z=r_k,$ where $r_k$ are roots of
 the following non-linear equation,
 \resizebox{1\columnwidth}{!}{
\begin{minipage}{1.28\columnwidth}
\begin{equation}\label{eq:delta_eq}
a^\top(Az+B)^{-1}a \,z + a^\top(Az+B)^{-1}b =0.
 \end{equation}
 \end{minipage}
}\\

\noindent Since $A,B$ are symmetric matrices from $\mathbb{R}^{2\times2}$, the matrix equation \eqref{eq:delta_eq} simplifies to
\resizebox{1\columnwidth}{!}{
\begin{minipage}{1.28\columnwidth}
\begin{equation}\label{eq:delta_eqsimple}
  a^\top  \text{adj}(A)   a\, z^2 + [a^\top \text{adj}(B)   a + a^\top \text{adj}(A)  b] z + a^\top  \text{adj}(B)  b=0, 
 \end{equation}
  \end{minipage}
}\\

 \noindent for such $z$  that $\det(Az+B)\neq0.$ Here, $\text{adj}(A):=D^\top\!\!AD$ is the adjugate matrix of $A.$
Depending on the value of the leading coefficient, Eq.~\eqref{eq:delta_eqsimple} is either a linear equation, if $a^\top  \text{adj}(A)  a=0,$ and has one root 
$$  r_1 =- \frac{a^\top \text{adj}(B)  a + a^\top \text{adj}(A)  b}{2 a^\top  \text{adj}(A)  a}, $$
or a quadratic equation  ($a^\top  \text{adj}(A)  a\neq0$) having at most two distinct roots:
$$
\begin{aligned}
r_1 = \frac{ - a^\top \text{adj}(B)  a - a^\top \text{adj}(A)  b  - \sqrt{d}} {2 a^\top  \text{adj}(A)  a },\\
r_2 = \frac{ - a^\top \text{adj}(B)  a - a^\top \text{adj}(A)  b  + \sqrt{d}} {2 a^\top  \text{adj}(A)  a },\\
\end{aligned}
$$
where $$d=[a^\top \text{adj}(B)  a + a^\top \text{adj}(A)  b]^2 - 4   a^\top  \text{adj}(A)  a  \,a^\top \text{adj}(B)  b.$$
Suppose, there is only one real root $r_1 \in [-1,1],$ which automatically implies that the other root either does not exist or is greater than 1 in its absolute value.
As being a convolution with the delta function, the integral  in \eqref{eq:integral} is simply an evaluation at a point,
$$ w_\infty(n+1) = \frac{1}{n}C_a[n (a  z+ b  + \frac{\mu_0}{n} )]C_b(n,z)\big|_{z=r1}.$$
After expanding $C_a(x),C_b(n,z)$ according to their definitions and some basic algebraic transformations the latter expression becomes
\resizebox{1\columnwidth}{!}{
\begin{minipage}{1.28\columnwidth}
\begin{equation}\label{eq:asymptote}
w_\infty(n+1) = L_0(L_1  n^{-\frac{3}{2} } +L_2 n^{-\frac{5}{2} } ) e^{ -\left(   E_1 n + E_0  +E_\text{-1}n^{-1}   \right) },\
\end{equation}
 \end{minipage}
}\\

\noindent and is exhaustively defined by definitions \eqref{eq:mu1mu2},\eqref{eq:AB} and the following list of constants:
\begin{equation}
\begin{aligned} 
L_0    &=2^{-3/2}\left(  \pi \det(S)\, a^\top S^{-1} a \right)^{-\frac{1}{2}},\\
L_1    &= C_1 \mu_0+  (r_1 a +b)^\top (C_2+C_2^\top) \mu_0  ,\\
L_2    &= \mu_0^\top C_2 \mu_0, \\
E_1    &=\frac{a^\top S^{-1}(a b^\top-b a^\top)S^{-1} b }{2a^\top S^{-1}a},\\
E_0    & =\frac{a^\top S^{-1}(a b^\top-b a^\top)S^{-1} \mu_0 }{a^\top S^{-1}a},\\
E_\text{-1} &=\frac{a^\top S^{-1}(a \mu_0^\top-\mu_0 a^\top)S^{-1} \mu_0 }{2a^\top S^{-1}a},\\
S      &=A r_1 + B. \\
\end{aligned}
\end{equation}
Note, that in the derivation of \eqref{eq:asymptote} the terms containing $n^{-0.5}$ cancel out.

If $E_1 \neq 0,$ the asymptote \eqref{eq:asymptote} decays exponentially fast,  and conversely,  $E_1 = 0$ is a sufficient and necessary condition
for the asymptote to decay as an algebraic function. The latter condition is equivalent to $$ a b^\top - b a^\top=0,$$ 
which after expansion according to the definitions \eqref{eq:mu1mu2},\eqref{eq:AB} simplifies to
$$2\mu \mu_{11}- \mu\mu_{02}   - \mu \mu_{20} + \mu_{02}   \mu_{20} - \mu_{11}^2=0.$$
This expression coincides with the definition of the critical point for the weak giant component \cite{kryven2016a}.

\subsection{Degenerate case of excess degree distribution}
Degenerated to the univariate case degree distributions, $u(k,l)=0,\;k>0,$ or $u(k,l)=0,\;l>0,$ present little interest as no connected components with size greater than 1 can be formed.
However, the asymptotic analysis for the case when one of the bivariate \emph{excess} distributions is degenerate, $u_\text{in}(k,l)=0,\;k>0$ or $u_\text{out}(k,l)=0,\;l>0$, requires a separate attention. 
Without loss of generality, suppose 
\begin{equation}\label{eq:degenerate_cond}
u_\text{in}(k,l)=0,\;k>0.
\end{equation}
Then, on one hand, covariance matrix $\Sigma_1$   is singular, and if  $\det (\Sigma_2) \neq0,$ the determinant $$\det (Az+B) = \frac{1}{2}\det [ (\Sigma_1-\Sigma_2 )z+\Sigma_1+\Sigma_2]=0$$  only  if $z=1.$ On another hand, $z=1$ is the only root of the quadratic equation \eqref{eq:delta_eqsimple} from the interval of interest, $z\in[-1,1].$ Consequently, Eq.~\eqref{eq:integral} fails to provide a valid description of the asymptote since $(Az+B)^{-1}$ does not exist at $z$, and one must seek an alternative route to perform the asymptotic analysis. 

Qualitatively, the condition \eqref{eq:degenerate_cond} means that there is at most one incoming edge per node. In view of the fact that the topology is locally tree-like, which is characteristic to finite components in configuration models, each finite component has exactly \emph{one} node with no incoming edges: the root node. 
Evidently, in this case, the asymmetry of the edges forces the connected components to be globally asymmetric as well: there is exactly one node per component with ingoing degree 0, and the whole component can be explored by starting at the root node and following exclusively outgoing edges.  We will now exploit the presence of such a global directionality in order to perform an asymptotic analysis for component sizes.

Let $w_0(n)$ denotes the probability that a component associated to the root node has size $n.$ It is $n$ times more likely to randomly select  any other node, then the root from a given component. Therefore $$w(n) =  \frac{1}{C} n w_0(n),$$ where  the normalisation constant $C$ is the expected component size. 
The condition on $u_\text{in},$ as given in Eq.~\eqref{eq:degenerate_cond}, can be rewritten as a condition on $u(k,l),$ that is  $u(k,l)=0,\;k>1.$ Let us introduce some auxiliary notation:
$$\mu_{0}'  = \sum\limits_{l=0}^\infty  u( 0, l ),\;\mu_{1}' = \sum\limits_{l=0}^\infty  u( 1, l ),$$
$$ u_0(l)=\frac{u(0,l)}{\mu_{0}'},\;u_1(l)=\frac{u(1,l)}{\mu_{1}'},$$
$$
\mu_{0j}' =\sum\limits_{l=0}^\infty l^j u_0(l),\; \mu_{1j}' =\sum\limits_{l=0}^\infty l^j u_1(l),
$$
where $j=0,1,2.$
We will go through a similar to Eq.~\eqref{eq:W1} derivation and construct a set of univariate equations for $w_0(n):$  
\begin{equation}\label{eq:W0}
\begin{aligned}
W_0(x)    & = x U_0[ W_{01}(x) ],\\
W_{01}(x) & = x U_1[ W_{01}(x) ],
\end{aligned}
\end{equation}
where $W_0(x),\; U_0(x),\; U_1(x)$ are generating functions for  respectively $w_0(n),\; u_0(l),$ and $u_1(l).$
By solving \eqref{eq:W0} for  $C=W_0(1)$ one obtains the expected component size 
$$C= 1 + \frac{\mu_{01}'}{1 - \mu_{11}'}.$$ 
%\frac{1}{\mu_0'}
Furthermore,  applying the Lagrangian inversion to Eq.~\eqref{eq:W0} gives the formal solution
$$
w_0(n)=\frac{1}{n-1}[k u_0(k)*u_1^{*n-1}(k)](n-2),
$$
which leads to the following asymptote for large $n$:
\begin{equation}\label{eq:degenerate_eq}
w_{\infty}(n) = L_{0}  n^{-\frac{1}{2}}  e^{ -E_0 - n E_1  },
\end{equation}
where
$$
\begin{aligned}
L_0 &=  \frac{ \mu_{01}'  (\mu_{11}' - 1 ) } { ( \mu_{11}' - \mu_{01}' -1)  \sqrt{ 2  \pi (\mu_{12}'-\mu_{11}'^2 ) }};\\
E_0 &= \frac{ (  \mu_{11}' -1 ) ( \mu_{01}' + \mu_{02}' -   \mu_{01}'  \mu_{11}' )  }  { \mu_{01}'  (  \mu_{12}'-\mu_{11}'^2  )};\\
E_1 &= \frac{( \mu_{11}' -1)^2}{2 (\mu_{12}'-\mu_{11}'^2  )}.
\end{aligned}
$$
In contrast to the non-degenerate asymptote \eqref{eq:asymptote}, that has the leading exponent $-\frac{3}{2}$, the leading exponent in the degenerate asymptote \eqref{eq:degenerate_eq} is $-\frac{1}{2}.$ 
That said, a pure algebraic asymptote $n^{-\frac{1}{2}}$ cannot be observed under the condition of finite moments $\mu_{03}', \mu_{12}'$. 
Indeed, $E_1\to0$ also implies that $\mu_{11}'-1\to 0,$ and consequently, pre-factor $L_0$ vanishes as well.

\section{Multiplex networks}
\subsection{General case of an arbitrary number of layers}
\begin{figure}[h]%%%%%%%%%
\begin{center}
\includegraphics[width=0.7\columnwidth]{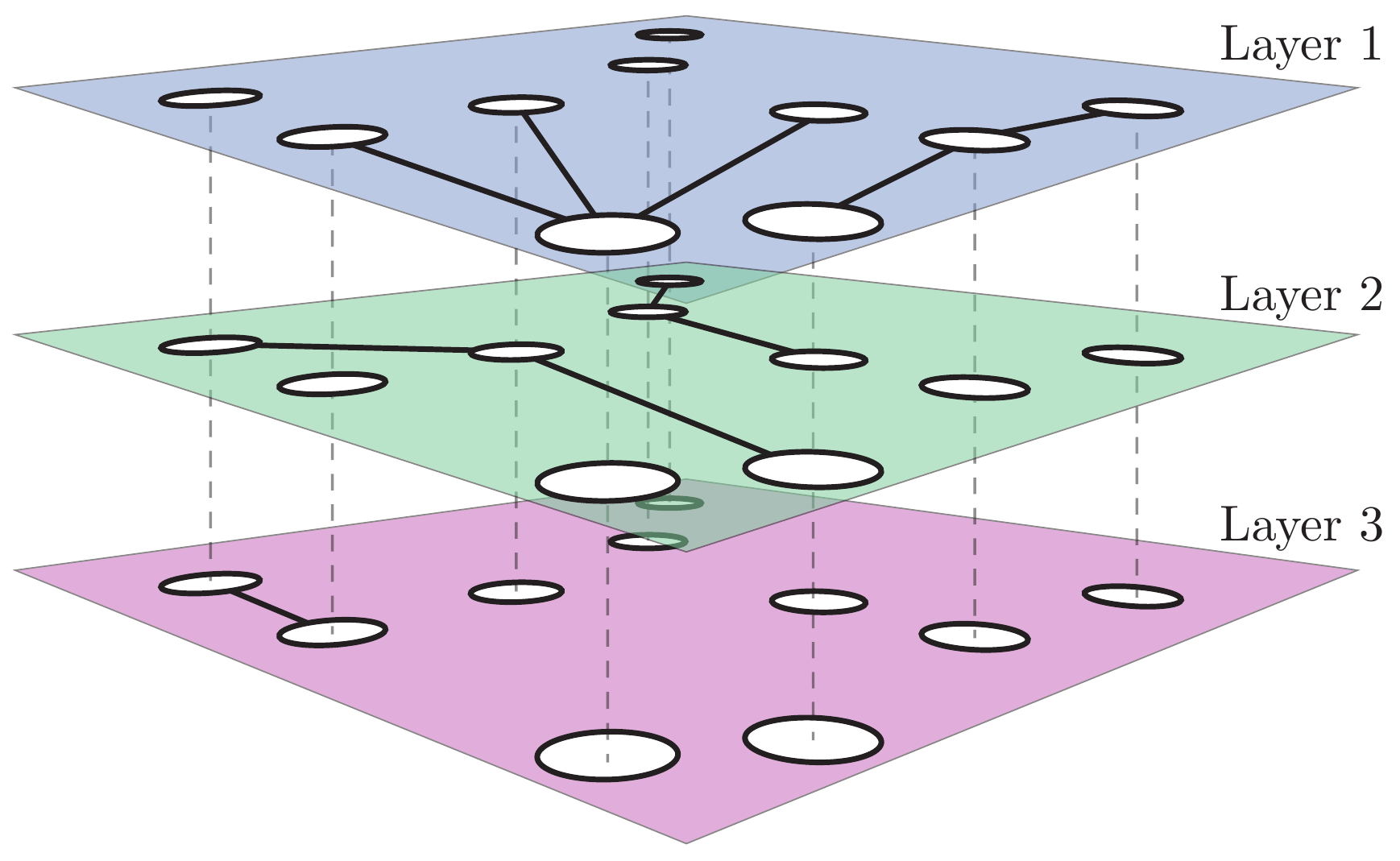}
\caption{An example of multiplex network with three layers. Each edge belongs to only one layer, whereas each node has one 'copy' in each layer. Even though, in each separate layer the nodes can be partitioned into different sets of connected components, the network is fully connected in the weak sence, when the connectivity information from all layers is combined. }
\label{fig:layers}
\end{center}
\end{figure}
This section considers the multiplex configuration model: a generalisation of the configuration model in which undirected edges are partitioned into subsets commonly referred to as types, layers, or colors \cite{kivela2014}.
In multiplex networks each edge belongs to one of many layers. Figure~\ref{fig:layers} illustrates an instance of a three-layer multiplex network with 10 nodes.
There are multiple ways to define a path in multiplex networks. A multilayer path, or simply a path in this section of the paper, is a path that combines edges from arbitrary layers. This definition of a path gives rise to a definition of multilayer connected components as sets of nodes connected together with the path.

Suppose, each edge belongs to a layer from $\Omega=\{1,\dots,N\}.$ We update the definition of the degree distribution to be a multivariate function $u(k_1,\dots,k_N),\; k_i\in \mathbb{N}_0$ that denotes probability of randomly choosing a node with  $k_i$ adjacent edges of from layer $i.$ 
As before, the degree distribution is normalised  $\sum\limits_{k_1,\dots,k_N} u(k_1,\dots,k_n)=1.$ 
The excess degree distribution associated to the $i^\text{th}$ layer is denoted with $u_i(k_1, \dots ,k_N) := \frac{ k_i + 1 }{\mu_i}u_i(k_1,\dots,k_i+1,\dots,k_N),$ 
where $\mu_i=\sum\limits_{k_1,\dots,k_N}(k_i+1) u_i(k_1,\dots,k_i+1,\dots,k_N)$ is the expected degree for $i^{\text{th}}$-layer edges.
Let $w(n)$ denotes the size distribution for multilayer components.  
By following similar considerations as in Section~\ref{eq:SectionDirrected}, one derives functional equations for the GF of the size-distribution, $W(n)$:
\begin{equation}\label{eq:ND1}
\begin{aligned}
W(x)&=x U[ W_1(x),\dots,W_N(x)],\\
W_1(x)&=x U_1[ W_1(x),\dots,W_N(x)],\\
\dots\\
W_N(x)&=x U_N[ W_1(x),\dots,W_N(x)],\\
\end{aligned}
\end{equation}
where the upper-case notation $U(x_1,\dots,x_N),$  $U_i(x_1,\dots,x_N), \;i=1,\dots,N$ is used to denote multivariate generating functions of the corresponding distributions.
The system of functional equations \eqref{eq:ND1} is a special case of \eqref{eq:Implicit_ND}, and thus can be solved by applying Lagrange-Good formula. 
Indeed, let $W_i(x)$ define the diagonals of $\mathbf A_i(\mathbf x)$, that is $\mathbf A_i(x,\dots,x) := W_i(x),\, |x|<1,\, x\in \mathbb{C}, $ for $ i=1,\dots,N.$  Additionally, let $ \mathbf R(\mathbf x)=[U_1(\mathbf x),\dots,U_1(\mathbf x)],$ and $F(\mathbf x) = U( \mathbf x),$ then the Lagrange-Good formula yields the expression for $a(k_1,\dots,k_N)$ that is generated by $A(x_1,\dots,x_N)=\frac{1}{x} W(x)$. The values for $w(n)$ can then be recovered using relation $w(n)=\sum\limits_{\substack{k_1+\dots+k_N=n-1\\k_i \geq 0}} a(k_1,\dots,k_N),$ and the complete equation for component-size distribution in the multilayered configuration network reads: for $n>1,$
\resizebox{1\columnwidth}{!}{
\begin{minipage}{1.28\columnwidth}
\begin{equation}\label{eq:multiplex}
w(n) =\sum_{\substack{k_1+\dots+k_N=n-1\\k_i \geq 0}} u(\mathbf k) *  \text{det}_*[D(\mathbf k)]* u_1(\mathbf k)^{*k_1}*\dots * u_N(\mathbf k)^{*k_N},
\end{equation}
  \end{minipage}
}\\
where 
$$D(\mathbf k)_{i,j}=\delta_{i,j} - [k_j u_i( \mathbf k ) ]* u_i( \mathbf k )^{*(-1)},\; i,j=1,\dots,N,$$
and $\text{det}_*[D(\mathbf k)]$ refers to the determinant of matrix $D$ computed with the multiplication replaced by the convolution: for example, 
$$det_*\begin{bmatrix}a(\mathbf k) & b(\mathbf k)\\ c(\mathbf k)& d(\mathbf k)\end{bmatrix}=a(\mathbf k)*d(\mathbf k)-b(\mathbf k )*c(\mathbf k).$$

\subsection{Two-layer multiplex network}
Suppose $N=2$, that is to say each edge belongs either to  Layer 1 or Layer 2. In this case, the degree distribution $d(k,l)$ is the probability of randomly selecting a node that bears $k$ edges in Layer 1 and $l$ edges in Layer 2. Where it leads to no confusion, we will reuse the notation from the previous section. For instance, the shorthand notations for moments, vectors of expected values and covariance matrices are as given in \eqref{eq:mom} and \eqref{eq:mu1mu2} respectively.  The total probability is normalised $\mu_{00}=1,$ but the expected numbers of edges in each layer need not be the same: $$\mu_{10}\neq \mu_{01}.$$  
The two dimensional version of \eqref{eq:ND1} reads,
\begin{equation}\label{eq:WW}
W(x) = x U\big[ W_{\text{1}}(x), W_{\text{2}}(x) \big],
\end{equation}
\begin{equation*}
\tag{\ref{eq:WW}$a$}
\begin{aligned}
&W_{\text{1}}(x) = x U_{\text{1}}\big[ W_{\text{1}}(x), W_{\text{2}}(x) \big],\\
&W_{\text{2}}(x)  = x U_{\text{2}} \big[ W_{\text{1}}(x), W_{\text{2}}(x) \big],
\end{aligned}
\end{equation*}
where $U(x,y),U_1(x,y),U_2(x,y)$ denote the corresponding generating functions for degree and excess distributions, and $W(x)$ is the generating function for the size distribution of two-layer connected components. The only structural difference between the equation for directed networks Eq.~\eqref{eq:W1} and the equation for two-layered network \eqref{eq:WW} is the order of arguments in the degree distribution GFs, which indicates presence or absence of symmetric edges (compare $W_{in}(x)=x U_{in}[W_{out}(x),W_{in}(x)]$ against $W_1(x) =x U_{1}[W_{1}(x),W_{2}(x)]$).
By setting $N=2$ in \eqref{eq:multiplex} one obtains the formal solution to \eqref{eq:WW}, 
\begin{equation}\label{eq:w_2layer}
w(n)= \sum\limits_{i=0}^{n-1} a(i,n-i-1),\; n>1,
\end{equation}
where for $i,j\geq 0$,
\begin{equation}
\label{eq:ajkL2}
a(i,j)= u(k,l)* u_{\text{1}}(k,l)^{*(i-1)}* u_{\text{2}}(k,l)^{*(j-1) } * d(k,l)\Big|_{\substack{k=i\\l=j}},
\end{equation}
and
\begin{multline}
\label{eq:DjkL2}
d(k,l)= [ u_{\text{1}}(k,l)- k u_{\text{1}}(k,l)] *  [ u_{\text{2}}(k,l) - l u_{\text{2}}(k,l) ] \\- l u_{\text{1}}(k,l) * k u_{\text{2}}(k,l).  
\end{multline}
We will now see how the asymptotic theory from Section~\ref{sec:asdirect} can be recast to fit the case of the two-layer multiplex networks.

\subsection{Asymptotic analysis for a bilayer network}
Let $\mu_1, \mu_2$ denote expected values and $\Sigma_1,\Sigma_2$ covariance matrices of $\frac{k}{\mu_{10}}u(k,l)$ and $\frac{l}{\mu_{01}}u(k,l),$ as given in definition \eqref{eq:mu1mu2}.
The characteristic function for the right hand side of Eq.~\eqref{eq:ajkL2} reads:
\resizebox{1\columnwidth}{!}{
\begin{minipage}{1.28\columnwidth}
\begin{equation}
\begin{aligned}
\label{eq:long_ch2}
\phi_a&(\omega_1,\omega_2)=  e^{-\iu (i \omega_1+ j\omega_2)} \phi(\omega_1, \omega_2)\\
\times& \left( 
\left[- \frac{\iu}{\mu_1}\frac{\partial}{\partial \omega_1}\phi(\omega_1, \omega_2)\right]^i
\left[-\frac{\iu}{\mu_2} \frac{\partial}{\partial \omega_2}\phi(\omega_1, \omega_2)\right]^j
\right.\\
&+
\frac{\iu}{j}  
\frac{\partial }{\partial \omega_2}\left[- \frac{\iu}{\mu_1} \frac{\partial}{\partial \omega_1}\phi(\omega_1, \omega_2)\right]^i
\left[- \frac{\iu}{\mu_2} \frac{\partial}{\partial \omega_2}\phi(\omega_1, \omega_2)\right]^j\\
&+
\frac{\iu}{i}  
\left[- \frac{\iu}{\mu_1} \frac{\partial}{\partial \omega_1}\phi(\omega_1, \omega_2)\right]^i
\frac{\partial }{\partial \omega_1}\left[ - \frac{\iu}{\mu_2}\frac{\partial}{\partial \omega_2}\phi(\omega_1, \omega_2)\right]^j\\
&-
\frac{1}{i j}  
\frac{\partial }{\partial \omega_2}\left[- \frac{\iu}{\mu_1} \frac{\partial}{\partial \omega_1}\phi(\omega_1, \omega_2)\right]^i
\frac{\partial }{\partial \omega_1}\left[- \frac{\iu}{\mu_2} \frac{\partial}{\partial \omega_2}\phi(\omega_1, \omega_2)\right]^j\\
&+
\left.
\frac{1}{i j} 
\frac{\partial }{\partial \omega_1}\left[- \frac{\iu}{\mu_1} \frac{\partial}{\partial \omega_1}\phi(\omega_1, \omega_2)\right]^i
\frac{\partial }{\partial \omega_2}\left[- \frac{\iu}{\mu_2} \frac{\partial}{\partial \omega_2}\phi(\omega_1, \omega_2)\right]^j\right).
\end{aligned}
\end{equation}
In the large $n$ limit, the latter approaches
\begin{equation}\label{eq:a_inf2}
a_{\infty}(i,j) =C(x) \frac{ e^{ -\frac{1}{2} x^\top\Sigma^{-1} x }}{2 \pi  \sqrt{\det(\Sigma)}},\\
\end{equation}
where $x = 2(i,j) - ( i\mu_1^\top + j\mu_2^\top)+\mu_0$ and
\begin{equation}\label{eq:C0C1C2}
\begin{aligned}
C(x) = &C_0+C_1  \frac{x}{n} + \frac{x^\top}{n} C_2\frac{x}{n},\\
C_0    = & 4 - 2 (I_1 \mu_1 + I_2 \mu_2) - \mu_1 D \mu_2,\\
C_1    = & [\mu_2^\top D \Sigma_1 - \mu_1^\top D \Sigma_2 - 2 ( I_1^\top \Sigma_1 + I_2^\top \Sigma_2  )]\left(\frac{1}{n}\Sigma\right)^{-1},\\
C_2    = &-\left(\frac{1}{n}\Sigma\right)^{-1} \Sigma_1 D \Sigma_2 \left(\frac{1}{n}\Sigma\right)^{-1} , \\
\Sigma=&i\Sigma_1+j\Sigma_2.
\end{aligned}
\end{equation}
 \end{minipage}
}\\

\noindent By applying the change of variables \eqref{eq:new_variables}, one obtains $ x = n ( a  z+ b  + \frac{\mu_0}{n} )$ and $\Sigma = n ( A z + B )$
where
\begin{equation}
\begin{aligned}\label{eq:abAB}
a &=  I_1-I_2+\frac{\mu_1 - \mu_2}2, \\
b &=  1 - \frac{\mu_1 + \mu_2}2,\\
A&=\frac{\Sigma_1- \Sigma_2}{2},\\
B&=\frac{\Sigma_1 + \Sigma_2}{2}.
\end{aligned}
\end{equation}
Now, coefficients $C_0,C_1,C_2$ and $\frac{1}{n}\Sigma = A z + B$ became independent of $n,$ and Eq.~\eqref{eq:a_inf2} is identical to \eqref{eq:a_inf} up to the definitions of constants $a,b,A,B,C_0,C_1,C_2,$ that are given above. Therefore, we can readily use the asymptote \eqref{eq:asymptote} also in the case of a two-layer networks. It is enough to redefine the constants according to definitions \eqref{eq:C0C1C2} and \eqref{eq:abAB} and take $z=r_1$, where $r_1\in[-1,1]$ denotes the root of Eq.~\eqref{eq:delta_eqsimple}.  As before, condition $ a b^\top - b a^\top=0$ indicates the emergence of the algebraic decay $n^{-\frac{3}{2}}$ in the sizes of connected components. When the latter equality is expanded according to definitions \eqref{eq:abAB} and \eqref{eq:mu1mu2}, one obtains the criterion in terms of degree distribution moments: 
\begin{equation}\label{eq:criterion_two_layers}
G(u):=\mu_{11}^2 -  (\mu_{20}-2 \mu_{10})(\mu_{02} - 2 \mu_{01}  ) = 0.
\end{equation}
As in the case of directed networks, the degenerate excess degree distribution, $u_1(k,l)=0,\;k>0,$ renders the asymptotic analysis not applicable. Nevertheless, the degenerate case is equivalent to a one-layer network with coupled nodes that has a univariate degree  distribution $d(l)=d(0,l)+\frac{1}{2}d(1,l),\; l=0,\dots$ The asymptotic theory for mono-layer components has been covered in Ref.~\cite{kryven2017a}, and, unlike in the case of directed networks, no new asymptotic modes emerge when the excess distribution is degenerate.

\subsection{Criticality in two-layer multiplex networks}\label{sec:criticality}
When a configuration network is a two-layered, one may speak of a connected component contained within a specific layer: a path comprise solely of edges from one layer, or a multilayer (weak) connected component that emerges from a combination of two layers: both types of edges may appear in the path.
No matter what type of connected components is considered, the asymptote of the component size distribution exhibits either exponential or algebraic decays.

When focusing on single-layer connected components, for instance, in Layer 1, condition $\mu_{20}-2\mu_{10}=0$ signifies the critical regime of the corresponding size distribution. Furthermore, a giant component exists within this layer iff $\mu_{20}-2\mu_{10}>0$.
Existence of a giant component within a single layer is a strong condition: it automatically implies existence of the weak two-layer giant component. More importantly, the two-layer giant component can also exist even when there are no single-layer giant components.

When two-layered connected components are considered, criterion \eqref{eq:criterion_two_layers} gives the condition for the algebraic decay in the component-size distribution. It is important to note that one should consider this inequality only together with the validity conditions of the asymptotic theory: $\mu_{ij}<\infty,\;i+j \leq 3 $ and existence of the root of Eq.~\eqref{eq:delta_eq}, $|r_1|\leq1$. For instance, unlike in the case of a single layer network, one cannot associate existence of the two-layer giant component solely with the sign of the left hand side of Eq.~\eqref{eq:criterion_two_layers}.  For a simple contra example, set $\mu_{11}=0.$ Then, the left hand side of Eq.~\eqref{eq:criterion_two_layers} is positive iff one layer contains a giant component and the other -- does not. When both layers contain a giant component simultaneously (or both layers contain no giant component), the sign is negative.  

Now, let us consider a critical degree distribution $u_c(k,l)$ such that $G(u_c)=0$ and $\mu_{11}>0.$ Assume, there are no single-layer giant components, that is to say $2\mu_{10}-\mu_{20}$ and $2\mu_{01}-\mu_{02}$ are positive quantities. The upper bounds on these quantities can be obtained from the Cauchy-Schwarz inequality, $\mu_{11}^2 \leq \mu_{20}\mu_{02}.$ The latter, when combined with condition $G(u_c)=0$ yields 
\begin{equation}\label{eq:single_layer_bounds}
0<2\mu_{10}-\mu_{20}\leq \frac{\mu_{10} \mu_{02} }{\mu_{01}},\;
0<2\mu_{01}-\mu_{02}\leq \frac{\mu_{01} \mu_{20} }{\mu_{10}}.
\end{equation}
Since there are no isolated nodes, $u(k,l)=0,\;k,l=0,$ the sum of expected numbers of edges is bounded below with
\begin{equation}\label{eq:m10m01sum}  
\mu_{10}+\mu_{01} \geq 1.
\end{equation}
Additionally, one obtains the following bounds from the monotonicity of the moments,
\begin{equation}\label{eq:monotonicity}
\mu_{20} < \mu_{10}^2,\; 
\mu_{02} < \mu_{01}^2. 
\end{equation}
Let us perturb expected number of edges $\mu_{10}$ by uniformly adding (or removing) a small number of edges $\text{d}\alpha$ in the first layer. Due to this perturbation, the degree distribution variates as $\text{d}u(k,l) = [u(k-1,l)-u(k,l)]\text{d}\alpha.$ The perturbation conserves the total probability, $\sum\limits_{k,l=0}^{\infty}\text{d} u(k,l) = 0,$ whereas the expected number of edges, indeed, variates as $\text{d}\mu_{10} = \sum\limits_{k,l=0}^{\infty} k [u(k-1,l)-u(k,l)] \,\text{d}\alpha =  \sum\limits_{k,l=0}^{\infty}[(k+1) u(k,l) - k u(k,l)]\,\text{d}\alpha = \text{d}\alpha$.  After expanding variations $\text{d}\mu_{11}=\mu_{01}\text{d}\alpha$ and $\text{d}\mu_{20}=(2\mu_{10}+1)\text{d}\alpha$ in a similar fashion,  we write the G\^ateaux derivative, 
\begin{multline*}
\frac{\text{d}}{d\alpha} G(u_c)=\lim\limits_{\text{d}\alpha\to 0}\frac{G(u_c + \text{d}u )-G(u_c)}{\text{d}\alpha} =  \\(\mu_{11}+\mu_{01}\text{d}\alpha)^2 - [\mu_{20}+(2\mu_{10}+1)\text{d}\alpha-2 (\mu_{10}+\text{d}\alpha)](\mu_{02} - 2 \mu_{01}  )  =\\
(2 \mu_{01} - \mu_{02}) (2 \mu_{10} - 1) + 2 \mu_{01} \mu_{11}. 
\end{multline*}
We will now show that $\frac{\text{d}}{d\alpha} G(u_c)>0$ by considering two cases.  Firstly, let $2 \mu_{10} - 1\geq 0,$ then  $(2 \mu_{01} - \mu_{02})>0$ according to \eqref{eq:single_layer_bounds}, and consequently, $\frac{\text{d}}{d\alpha} G(u_c)>0.$
Secondly, let us assume the opposite is true, $ 0<\mu_{10} < \frac{1}{2}:$ by expressing $\mu_{11}$ from \eqref{eq:criterion_two_layers} and plugging it into 
$\frac{\text{d}}{d\alpha} G(u_c)>0$ one obtains
\begin{equation}\label{eq:temp}
\mu_{01} (2 \mu_{10} - \mu_{20}) - \frac{1}{\mu_{01} }(2 \mu_{01} - \mu_{02}) (1 - 2 \mu_{10})^2>0.
\end{equation}
Combining lower bound $\mu_{01}\geq1-\mu_{10}=\frac{1}{2}$ (as follows from \eqref{eq:m10m01sum}) and the upper bound on $\mu_{20}$ (as given in \eqref{eq:monotonicity}), the first term in \eqref{eq:temp} is bounded below with $\mu_{01} (2 \mu_{10} - \mu_{20}) \geq \frac{1}{2} (2 \mu_{10} - \mu_{10}^2). $
The lower bound for the second term of \eqref{eq:temp} follows from sequentially applying \eqref{eq:single_layer_bounds} and \eqref{eq:monotonicity}:
\begin{multline}
- \frac{1}{\mu_{01} }(2 \mu_{01} - \mu_{02}) (1 - 2 \mu_{10})^2  \geq \\ -\frac{ \mu_{20} }{\mu_{10}}(1 - 2 \mu_{10})^2 \geq \\ -\mu_{10} (1 - 2 \mu_{10})^2.
\end{multline} So that
\begin{multline}
\mu_{01} (2 \mu_{10} - \mu_{20}) - \frac{1}{\mu_{01} }(2 \mu_{01} - \mu_{02}) (1 - 2 \mu_{10})^2 \geq \\  \frac{1}{2} (2 \mu_{10} - \mu_{10}^2)-\mu_{10} (1 - 2 \mu_{10})^2= \frac{1}{2} (7 - 8 \mu_{10}) \mu_{10}^2> \frac{3}{2} \mu_{10}^2>0.
\end{multline}
The fact that $\frac{\text{d}}{d\alpha} G(u_c)>0$ means that perturbing the configuration network at the critical regime, $G(u)=0,$ by a uniform addition of new edges forces the value of $G(u)$ to become positive.  The opposite is also true: uniform removal of existing edges at the critical regime forces values of $G(u)$ to become negative. Similar derivation holds for the perturbation in the second layer.

Finally, suppose one modifies $u(k,l)$ in such a manner that the expected numbers of edges, $\mu_{10},\mu_{01}$ remain constant whereas the second moments alter. Such a perturbation of the degree distribution causes rewiring of the network but keeps the expected numbers of edges in each layer the same. Consider a function $f( \mu_{20}, \mu_{02} ):=\mu_{11}^2 -  (\mu_{20}-2 \mu_{10})(\mu_{02} - 2 \mu_{01}  )=G(u_c),$  as follows from the lower bounds \eqref{eq:single_layer_bounds}, both components of  the gradient vector,
$$ \nabla f( \mu_{20}, \mu_{02} ) = (2 \mu_{01} - \mu_{02}, 2 \mu_{10} - \mu_{20}),$$
are positive. This fact confirms that rewiring that moves edges within a single layer toward the nodes with higher degree forces the value of $G(u)$ to become positive. The total action of the simultaneous rewiring in two layers is defined by  
$\text{sign}[ (2 \mu_{01} - \mu_{02}, 2 \mu_{10} - \mu_{20})^\top (\partial \mu_{20},\partial \mu_{02})].$

According to the asymptote \eqref{eq:asymptote}, if $G(u)=0$, the size distribution decays algebraically with exponent $-\frac{3}{2},$ and therefore, expected component size diverges. On one hand, perturbations of the network by inflection with new edges or a rewiring that moves edges to nodes with larger degree does not reduce the size of the largest component, on the other hand, after such a perturbation $G(u)\neq 0$: the component-size distribution switches to the exponential decay and features finite expected component size. The deficit in expected component size, which due to the nature of the perturbation could have only increased if all connected  components were finite, is attributed to the emergence of the giant two-layer component. 
\begin{figure}
\begin{center}
\includegraphics[width=\columnwidth]{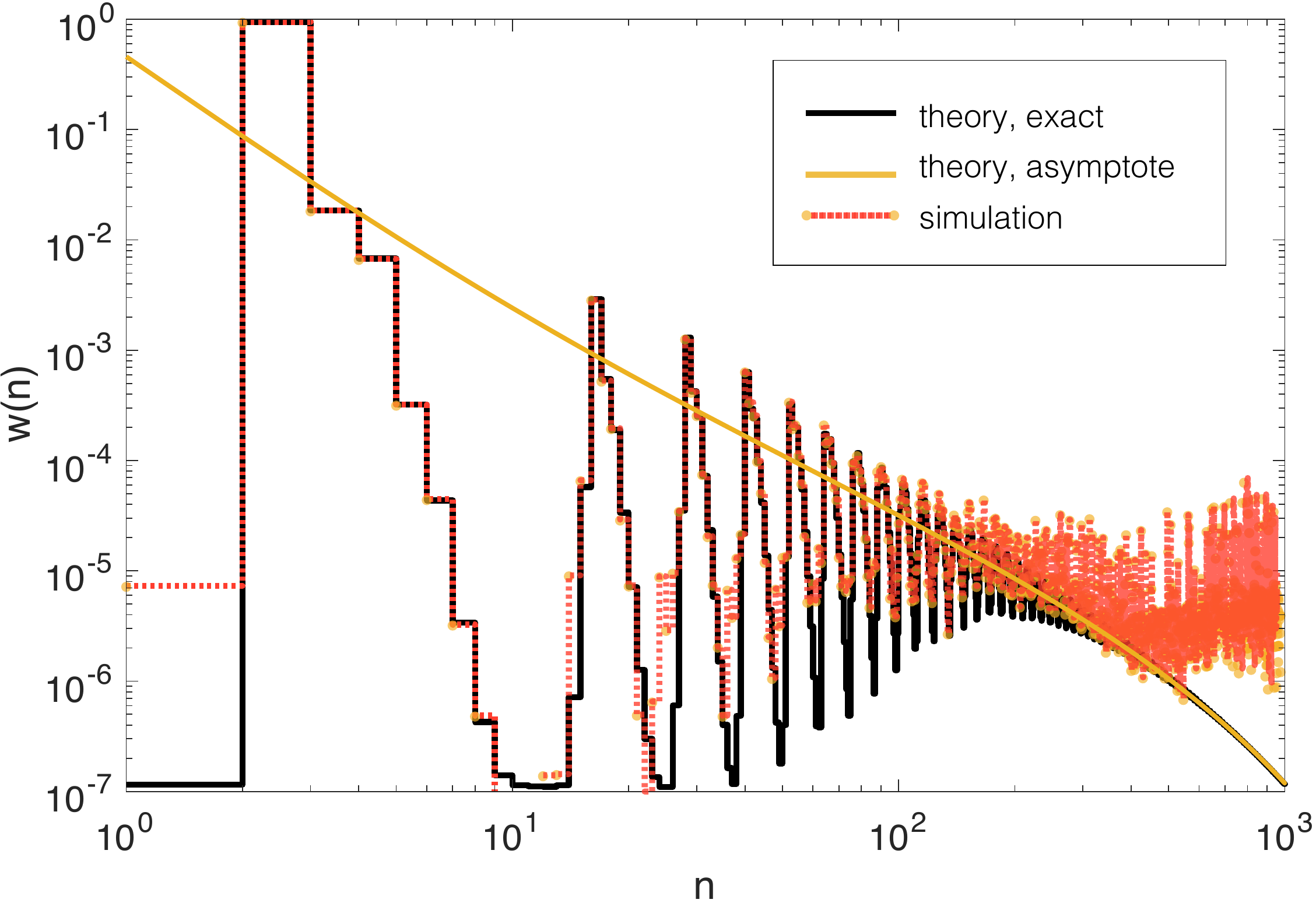}
\caption{(\emph{color online})
 An oscillatory example of the  size distribution of connected components in the two-layer configuration model as predicted by the analytical expression \eqref{eq:w_2layer} (\emph{solid black line}) is compared against the data obtained from simulations  (\emph{scattered points} linked with a \emph{dashed line} that indicates the trend).
Low probabilities are naturally underrepresented in simulated data due to a limited size of the Monte Carlo sample. 
The theory, as given by Eqs. \eqref{eq:asymptote} and \eqref{eq:C0C1C2}, predicts the asymptote with a transient slope $-\frac{3}{2}$ (\emph{solid yellow line}).
}
\label{fig:two_layer}
\end{center}
\end{figure}
\begin{figure}
\begin{center}
\includegraphics[width=\columnwidth]{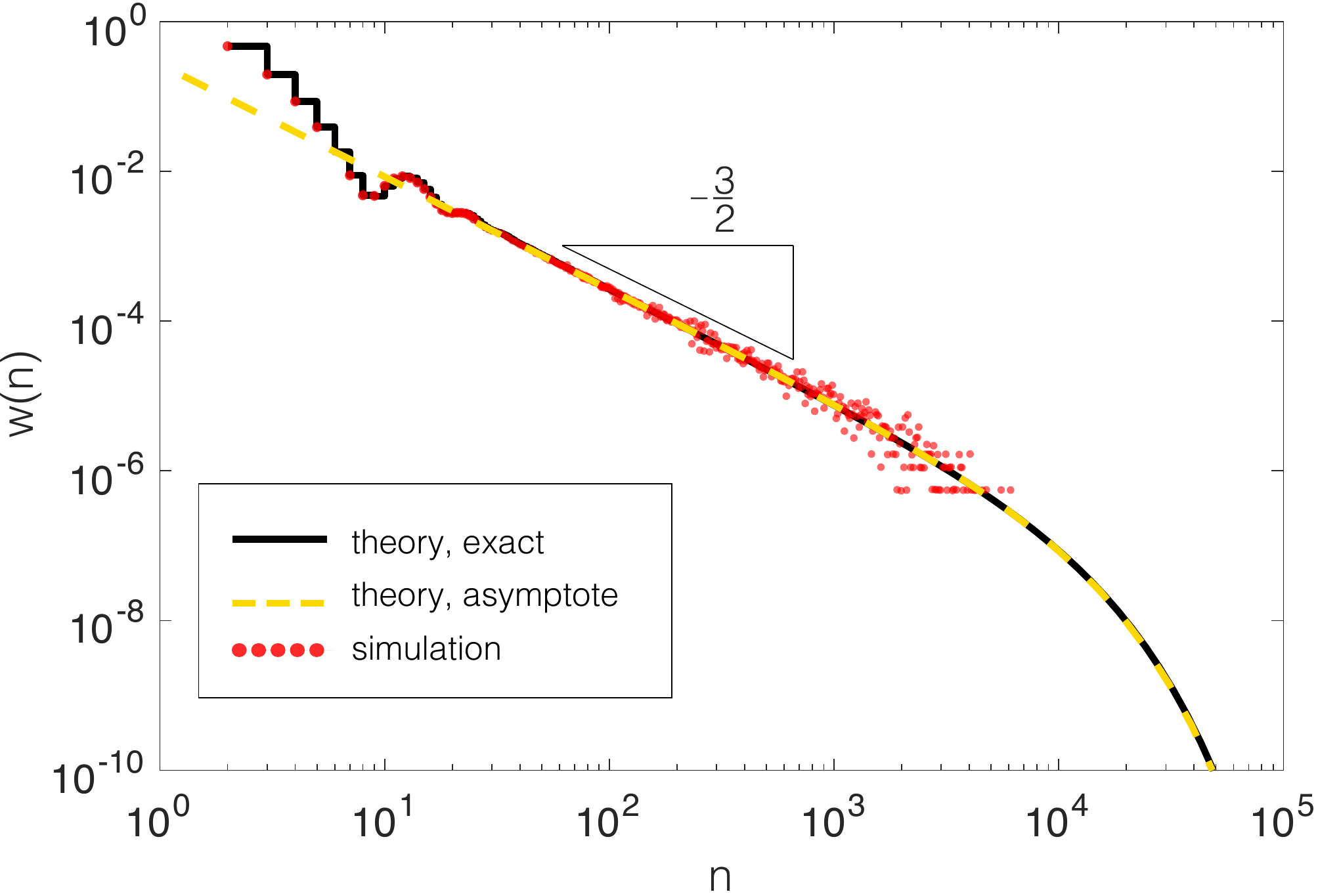}
\caption{(\emph{color online}) An example of the size distribution of weakly connected components in a directed configuration model as predicted by the analytical expression \eqref{eq:w2d} (\emph{solid line}) is compared against the  data obtained from simulations  (\emph{scattered points}).
The theory, as given by Eq. \eqref{eq:asymptote}, predicts the asymptote with a transient slope $-\frac{3}{2}$ (\emph{dashed line}).}
\label{fig:32}
\end{center}
\end{figure}
\begin{figure}
\begin{center}
\includegraphics[width=\columnwidth]{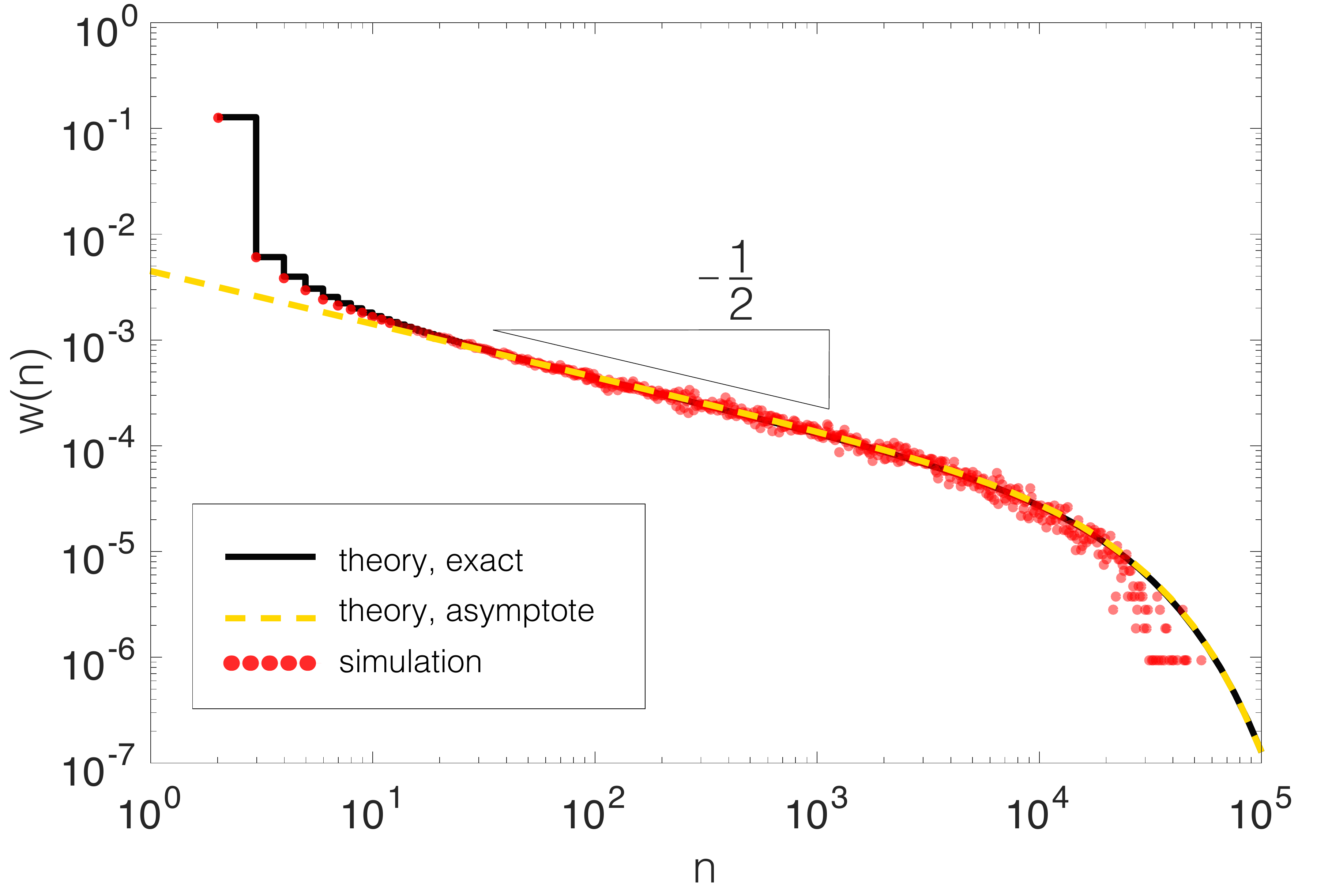}
\caption{(\emph{color online}) An example of the size distribution of weakly connected components in a directed configuration model as predicted by the analytical expression \eqref{eq:w2d} (\emph{solid line}) is compared against the  data obtained from simulations  (\emph{scattered points}).
The theory, as given by Eq. \eqref{eq:degenerate_eq}, predicts the asymptote with a transient slope $-\frac{1}{2}$ (\emph{dashed line}).}
\label{fig:12}
\end{center}
\end{figure}

\section{Discussion and conclusions}
The main results of this study are the formal expressions for the size distributions of connected components in directed and multiplex networks.
These expressions involve the convolution power and, in practice, can be  evaluated exactly with FFT algorithm in the cost of $O(n^2\log n),$ in the case of directed networks, and $O(n^{N}\log n ),$ in the case of multiplex networks with $1<N<n$ layers. These expressions are very general and do not rely upon any restrictions on the degree distribution itself. The supporting code is accessible at GitHub repository \cite{git:MultiDirNet}. 
Unlike the fixed point formulations \eqref{eq:W1},\eqref{eq:WW}, the formal expressions for component size distributions are tractable from asymptotic theory point of view. The asymptotic analysis for weak/multilayer connected components resulted in simple analytical expressions that, under certain conditions, feature a self-similar behaviour. 
The asymptotic theory, however, does rely on a few assumptions that to a certain extent limit the space of applicable degree distributions. First, we assume finiteness of partial moments, $\mu_{ij}<\infty,\; i+j \leq 3;$
second, we rely upon existence of the real root of Eq.~\eqref{eq:delta_eq} such that $|r_1| \leq 1.$
Finally, there is a practical restriction that arises if one aims to utilise the asymptote as an approximation for the size distribution itself: the best approximation accuracy is gained when the network is in the critical window
$|a r_1 +b| = \epsilon,$
where $\epsilon$ is infinitesimal.  

A few examples of size distributions of connected components and corresponding to them analytical asymptotes are given in Figs.~\ref{fig:two_layer},~\ref{fig:32} and \ref{fig:12}. Figure~\ref{fig:two_layer} compares the new theory against simulated data for the case of a two-layer network with the degree distribution  given  by 
$$
u(k,l) =  0.9782 e^{ - 5 [ ( k - 1)^2 + l^2 )  ]} +  0.002 e^{ - 10 [ ( k - 9 )^2 + ( l - 3 )^2]  } .
$$
This example was selected to demonstrates the possibility of an oscillatory behaviour arising in the size distribution of connected components.
It can be noted that the theory, as given by Eq.~\eqref{eq:w_2layer}, accurately  predicts the non-trivial oscillations present in the data.
 For large $n$, the theoretical predictions in this example converge to the asymptote, as given by Eq.~\eqref{eq:asymptote}.

Figure~\ref{fig:32} features the size distribution of connected components in a directed network featuring a non-degenerate degree distribution, 
$$
u(k,l) =   0.5167  e^{ - k^2 - l^2} + 0.0052 e^{ - 2.5[ ( k - 4 )^2 + ( l - 4 )^2 ] },
$$
whereas 
Fig.~\ref{fig:12} features the results obtained for a degenerate degree distribution:
$$
\begin{aligned}
u( 0, k ) &=     0.9073e^{ -2.266 k  },\;k\geq0,\\
u( 1, k ) &=     0.9073e^{ -0.7 k  },\;k\geq0,\\
u(l,k)    &= 0,\;l>1,k\geq0.
\end{aligned}
$$
As in the previous example, both Figs.~\ref{fig:32} and \ref{fig:12} compare the theoretical size distribution, as given by Eq.~\eqref{eq:w2d}, to the simulated data. In both figures,  the theoretical predictions and the data converge to the asymptotes for large $n$. 
In the case of the non-degenerate degree distribution, the asymptote features transient slope $-\frac{3}{2}$, as predicted by Eq. \eqref{eq:asymptote}.
However, in the case of the degenerate degree distribution the transient slope of the asymptote is $-\frac{1}{2}$, which is in accordance with Eq. \eqref{eq:degenerate_eq}. The latter observation is a surprising result. This is the first evidence that a configuration model with a light-tailed degree distribution may feature a distinct from $-\frac{3}{2}$ exponent. Importantly, in the multiplex configuration network with two layers such an anomaly is not present. When the degree distribution is light-tailed, both non-degenerate directed networks and two-layer networks feature leading exponent $-\frac{3}{2}$ in the critical regime, which is also the case in undirected networks.

\begin{figure}[H]
\begin{center}
a.\parbox{0.93\columnwidth}{\includegraphics[width=0.9\columnwidth]{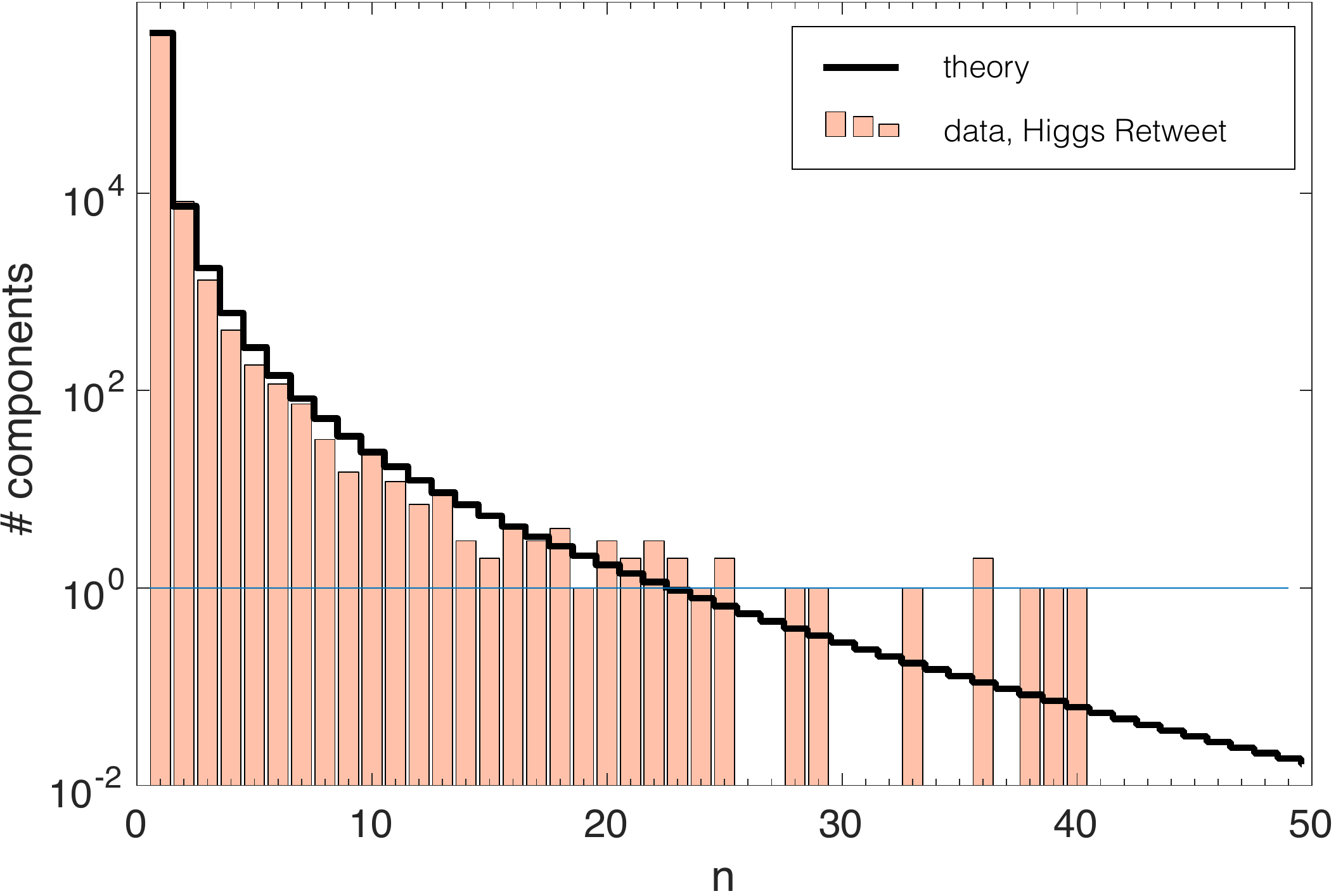}}
b.\parbox{0.93\columnwidth}{\includegraphics[width=0.9\columnwidth]{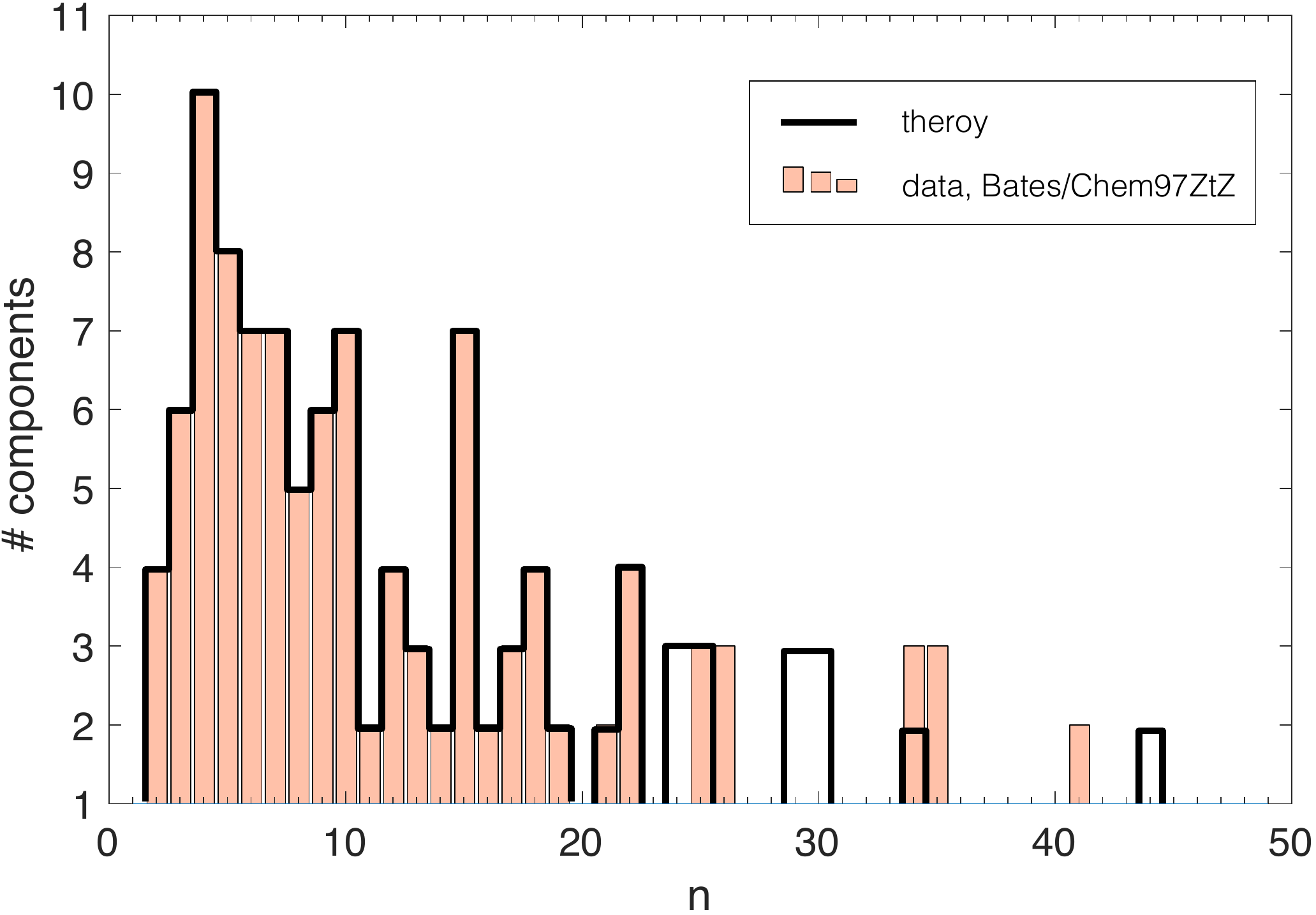}}
c.\parbox{0.93\columnwidth}{\includegraphics[width=0.9\columnwidth]{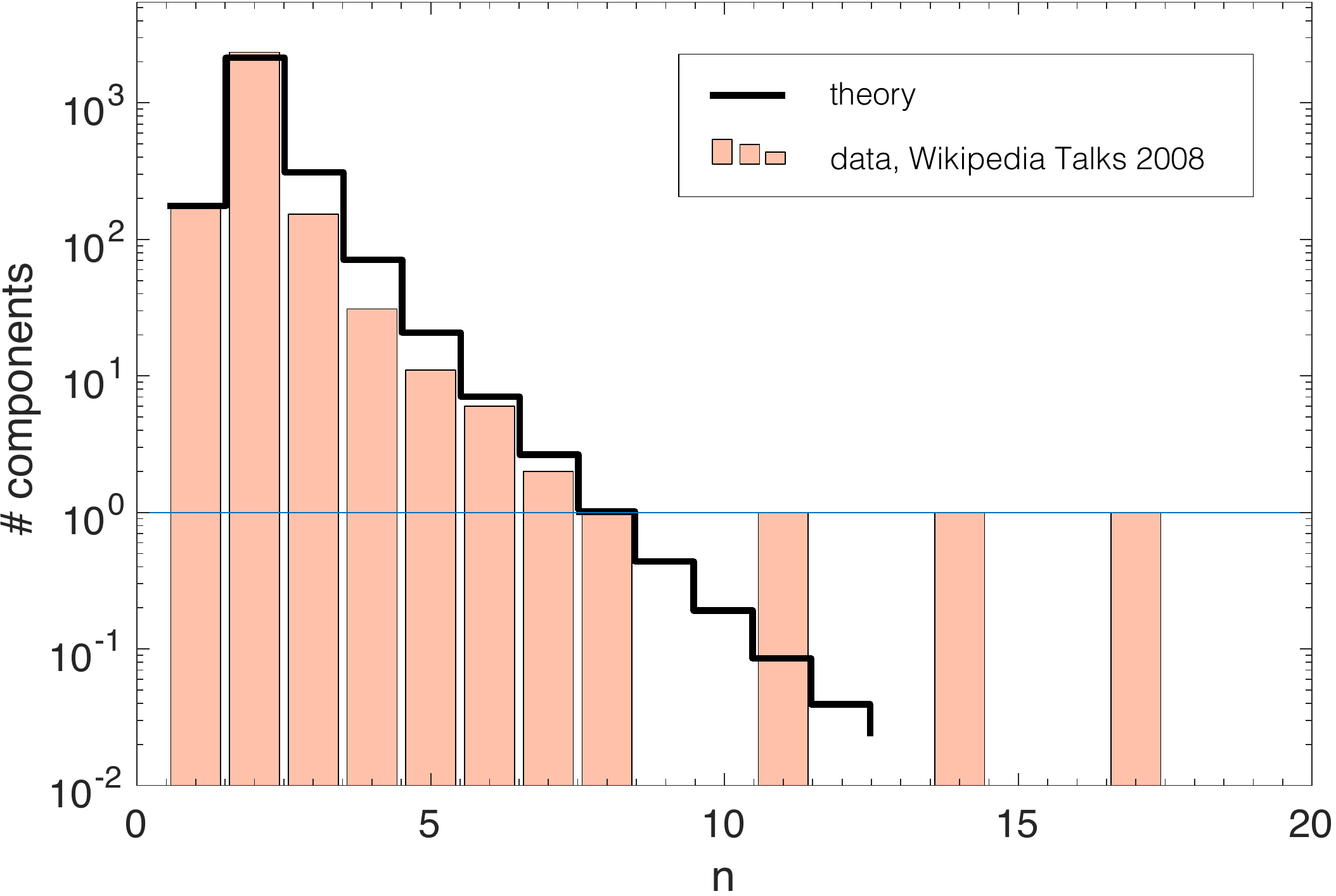}}
\caption{(\emph{color online}) Comparisons of respective theoretical size distributions of finite components against the empirical data. Three cases of directed networks containing $N$ nodes in total are considered:
\emph{a.} the network of retweets in the Higgs/Twitter dataset, $N=425\,008$ \cite{Manlio2013};  
\emph{b.} the graph of the sparse statistical matrix Chem97ZtZ, $N=2\,541$ \cite{davis2011};
\emph{c.} the network of communications on Wikipedia till January 2008, $N=2\,394\,385$ \cite{leskovec2010}.
}
\label{fig:empirical}
\end{center}
\end{figure}

A comparison of the theory agains a few examples of empirical data is given in Fig.~\ref{fig:empirical}. 
This figure presents normalised to the number of nodes theoretical size distributions of weakly connected components and compares them to empirical component-count distributions extracted from various datasets of directed networks.

In undirected, single-layer configuration networks, a heavy tail in the size distribution is observed when $\mu_2-2\mu_1=0$, where $\mu_2$ and $\mu_1$ are the moments of the univariate degree distribution. Furthermore, when the equality sign in this criterion is replaced by the inequality sign, $\mu_2-2\mu_1>0,$ one obtains the criterion for  giant component existence. Similar inequality criterion can be constructed for directed networks: also in this case the condition for a heavy tail in the size distribution relates to  the giant component existence \cite{kryven2016a}. However, this principle breaks down already in multiplex networks that consist of as few as two layers. The sign of the left hand side of the criticality condition \eqref{eq:criterion_two_layers} cannot be directly associated with existence of the giant two-layer component. Nevertheless, as was argued in Section \ref{sec:criticality}, if equality~\eqref{eq:criterion_two_layers} fails to hold due to a small perturbation in the expected numbers of edges $\mu_{10},\mu_{01}$ (or rewiring caused by an increase of second moments $\mu_{20},\mu_{02}$) at the critical regime, one can still associate the sign of the left hand side in Eq.~\eqref{eq:criterion_two_layers} with the giant component existence. This association is guaranteed  to be valid within the critical window. 

\begin{acknowledgments}
This work is part of the project number 639.071.511, which is financed by the Netherlands Organisation for Scientific Research (NWO) VENI.
\end{acknowledgments}

\bibliographystyle{apsrev4-1}
\bibliography{literature}

\end{document}